\documentclass[12pt, leqno]{article}
\usepackage{hyperref}
\usepackage{amsmath}
\usepackage{amsfonts}
\usepackage{amssymb}
\usepackage{graphicx}
\usepackage{color}
\usepackage[latin1]{inputenc}
\newtheorem{thm}{Theorem}[section]

\numberwithin{equation}{section}

\newcommand{\re}{\mathbb{R}}
\newcommand{\ren}{\mathbb{R}^n}
\newcommand{\Hn}{\mathbb{H}^n}
\newcommand{\Sn}{\mathbb{S}^n}

\newcommand{\ve}{\varepsilon}

\renewcommand{\dfrac}{\displaystyle\frac}
\definecolor{darkblue}{rgb}{0.05, .05, .65}
\definecolor{darkgreen}{rgb}{0.05, .70, .05}
\definecolor{darkred}{rgb}{0.8,0,0}
\def\qed{\unskip\kern 6pt \penalty 500
\raise -2pt\hbox{\vrule \vbox to8pt{\hrule width 6pt
\vfill\hrule}\vrule}\par}
\begin{document}
\title{\textbf{Fundamental solution and long time \\ behaviour of the Porous Medium \\
Equation  in Hyperbolic Space  }}
\author{Juan Luis V\'azquez\footnote{Dpto. de Matemáticas, Univ. Autónoma de Madrid, Spain.}
}
\date{ }
\maketitle

\begin{abstract}
We construct  the fundamental solution of the Porous Medium Equation posed in the hyperbolic space $\Hn$ and describe its asymptotic behaviour as $t\to\infty$. We also show that it describes the long time behaviour of integrable nonnegative solutions, and very accurately if the solutions are also radial and compactly supported. By radial we mean functions depending on the space variable only through the geodesic distance $r$ from a given point $O\in \Hn$. We show that the location of the free boundary of compactly supported solutions grows logarithmically for large times, in contrast with the well-known power-like growth of the PME in the Euclidean space. Very slow propagation at long distances is a feature of porous medium flow in hyperbolic space. We also present a non-uniqueness example for the Cauchy Problem based  on the  construction of an exact generalized traveling wave solution that originates from one point of the infinite horizon.
\end{abstract}

\section{\bf Introduction}

A number of degenerate diffusion equations exhibit the property of finite propagation whereby solutions with initial data localized in a region may expand the support in time, but only to a finite distance of the original support in each time interval. Such is the case of the porous medium equation (PME): $\partial_t u=\Delta u^m$, $m>1$, typically posed in the whole Euclidean space $\ren$, $n\ge 1$,  or in bounded subdomains thereof. Other examples are the Stefan Problem, the Hele-Shaw problem, or the $p$-Laplacian equation with $p>2$, cf \cite{Fr82}, \cite{VazICM}, among the many references.

Let us assume for simplicity that we deal with nonnegative solutions. The property of finite propagation gives rise in these models to the occurrence  of a free boundary $\Gamma$, which is the closed space-time set that separates the occupied or fluid region $\{u>0\}$   from the empty region $\{u=0\}$. The study of the location and regularity of free boundaries is a major topic in the research on nonlinear parabolic equations of degenerate type. It is usually proved that the free boundary is a nice hypersurface, except for some possible singular points.

A very representative situation  happens when we consider as initial data a Dirac mass, $u(x,0)=M\,\delta(x)$. When the linear heat equation (i.\,e., the case $m=1$ of the above porous medium example) is posed in the whole Euclidean space $\ren$, such data give rise to the Fundamental Solution, which has an explicit formula:
\begin{equation}\label{fs.he}
U(x,t)=M\,(4\pi t)^{-n/2}\,e^{-\frac{|x|^2}{4t}}\,,
\end{equation}
i.\,e.,  a Gaussian function with  variable variance, $\sigma^2(t)=2t$, which is the Brownian scaling. Note that in this equation nonnegative solutions are positive everywhere and there is no free boundary.

Free boundaries occur for the porous medium equation, $m>1$. If posed  for $x\in\ren$ and with the same type of initial data, we obtain instead of the above fundamental solutions  the famous family of Barenblatt solutions, first found and analyzed  by Zel'dovich, Kompanyeets and Barenblatt in Moscow around 1950   \cite{Bar1952, ZK1950}. These solutions have played a key role in the  theory of the equation, cf. \cite{Vbook}. It is a fortunate fact that such remarkable special solutions have a self-similar form,
\begin{equation}\label{Bar.euc.1}
U(x,t)=t^{-n\beta}F(x\,t^{-\beta})
\end{equation}
with exponent $\beta=1/(n(m-1)+2)$, and also that the profile is explicitly given by
\begin{equation}\label{Bar.euc.2}
F(\xi)=\left(C-k|\xi|^2\right)_+^{1/(m-1)}
\end{equation}
Here, $k=(m-1)\beta/2m$ is  fixed  while $C>0$ is an arbitrary constant that can be determined by fixing the mass of the solution $M=\int U(x,t)\,dx$ (which is invariant in time). The formula means in particular that the free boundary is the set
$$
\Gamma(U)=\{(x,t): t>0, |x|=r_0t^{\beta}\},
$$
where $ r_0^2=C/k=c(m,n)M^{2(m-1)\beta}$. The rate $|x|\sim t^{\beta}$ indicates the precise penetration strength of the porous medium flow in Euclidean space. Note that for $m>1$ we always have $\beta<1/2$; even more, $\beta\to 0$ as either $m$ or $n$ go to infinity, which means very weak propagation in such limit cases.

In this paper we ask the same kind of questions for the Porous Medium Equation posed in Hyperbolic space $\Hn$:
\begin{equation}\label{hpme}
\partial_t u=\Delta_g u^m, \ \quad m>1,
\end{equation}
that we label HPME for short;  here, $u=u(x,t)$, $x\in\Hn$, $t>0$, and $\Delta_g$ denotes the Laplace-Beltrami operator in the hyperbolic metric. Here are some typical questions: are there fundamental solutions? are they explicit or semi-explicit? can we calculate or estimate the height of the solution and the size of the support for $t>0$, at least for all large $t$? Is the fundamental solution representative of a wide class of other solutions? Are there any other important special solutions?

We will try to find answers to those questions below. In particular, we will see that the behaviour of the fundamental solutions for short and long  times is completely different, which is a remarkable difference with the Euclidean case described by (\ref{fs.he}) or (\ref{Bar.euc.1}). Our results demonstrate the gradual influence of curvature on the form of the fundamental solution in hyperbolic space: they start with a Barenblatt shape and evolve into a peculiar log-conical profile for very large times.


\subsection{\bf Organization of the paper and main results}

\noindent We start the study by recalling the formulas for the heat equation on hyperbolic space, as compared with the Euclidean case, Section \ref{sect.2}. In Section \ref{sect.3} we  introduce the porous medium equation in hyperbolic space and derive some of the basic estimates on the solutions  that will provide a rough idea of how solutions behave.  The main results for fundamental solutions are as follows

\begin{thm}\label{thm.main} Let $n\ge 3$. (i) Given an origin of radial coordinates $O$ there exists a unique radial and nonnegative weak solution $U(r,t)$ of the PME in $\Hn$ with initial data a unit delta function located at $O$. It is  bounded and continuous for $t>0$, monotone nonincreasing in $r$ for fixed $t> 0$, and it has compact support in a ball of geodesic radius $R(t)$ around $O$, that expands with time.

\noindent (ii) For small times it behaves like the Barenblatt solution of the Euclidean PME.

\noindent (iii) For $t\gg 1$ the behaviour of the fundamental solution is approximately given by the formula: \ $t\,U(r,t)^{m-1}\sim  a \left(\gamma \log t- r + b\right)_+$\,, where
\begin{equation}
 a=\frac{1}{m(n-1)}\,, \quad \gamma=\frac1{(m-1)(n-1)}\,,
\end{equation}
and \ $b(n,m)$ is fixed by the mass of the solution. More precisely, we have  convergence  along curves $r=\gamma\log (t)+\xi$ in the sense that
\begin{equation}\label{fundsol.hyper.1}
t\,U^{m-1}(\gamma\log(t)+\xi,t) \to a \left( b-\xi\right)_+
\end{equation}
as  $t\to\infty$ with fixed  $\xi$.  The convergence is uniform for $\xi\ge k$.

\noindent (iv) Moreover, there is a  free boundary $r=R(t)$ that grows for large times like $R(t)\approx \gamma \,\log (t)+b$. More precisely, \ $ e^{R(t)-b}\approx t^{\gamma}\,.$

\noindent (v) Finally, we have the long-time sup estimate:
\begin{equation}\label{supest}
\|U(\cdot,t)\|_\infty\sim \left(\log(t)/t\right)^{1/{m-1}}\,.
\end{equation}
\end{thm}

A word about notation: here and the sequel, $(f)_+=\max\{f,0\}$ means the positive part of $f$; the $\sim$ sign means that the ratio of the left-hand side to the right-hand side is bounded above and below by  positive constants. We will usually make precise the dependence of those constants in the paper. We will use the sign $\approx $ when the ratio tends to 1 in the corresponding limit. We will use the equivalent notations $u_t$ and $\partial_t u$ for the time derivative of $u$ depending on convenience, and similarly for the space derivatives.

Another precision:  the theorem refers to the fundamental solution with unit mass; a  fundamental solution with initial data $M\delta(x)$, $M>0$, can be obtained by simple scaling of the solution $U$ with unit mass according to the formula \ $U_M(r,t)=M\,U(r,M^{m-1}t)$. This is true both in the Euclidean and in the hyperbolic case.
In the latter case, the asymptotic formula (\ref{fundsol.hyper.1}) holds for $U_M$ with the same $a$ and $\gamma$, while $b$ is replaced by
\begin{equation}
b(M)=b+\gamma (m-1)\log M\,.
\end{equation}

Actually, in order to better understand the formulas in the main theorem, it is convenient to introduce the PME {\sl pressure}, a variable defined as $p=(m/(m-1))u^{m-1}$. This has been an important quantity in the porous medium theory in the Euclidean space since it allows for a common treatment of all exponents $m>1$; therefore, it has received considerable attention, cf. \cite{ArBk86, Vbook} for instance. The basic properties apply also in the hyperbolic space setting. In our concrete situation, the  results of Theorem \ref{thm.main} amount to say that the pressure of the fundamental solution satisfies the asymptotic estimate
\begin{equation}\label{press.formula}
P(r,t)\sim \frac{\gamma}{t}(\gamma \log(t)-r+b)_+\,.
\end{equation}

It is to be noted that our results do not apply exclusively to fundamental solutions. Generally, we deal with solutions with nonnegative integrable data, and very specially with radial solutions, by which we mean solutions depending only on the geodesic distance $r$ from a given point $O\in \Hn$ (and on time).

\medskip

\noindent $\bullet$ Going back to the outline of the paper, in Section
\ref{sect.tw} we prepare the way for the proof of the main results by first calculating an explicit solution (or more precisely, a family of solutions, see formula (\ref{gtw.mn})), that does not belong to the class of functions we are looking for (it is not integrable over the space) but will give us precise information of what to expect. The solution is interesting in itself as the following result shows.

\begin{thm}\label{thm1.2} There exists  nonnegative weak solution  of the HPME defined on the whole of $\Hn$ for all $t>0$, that has zero initial trace at $t=0$ and an expanding support bounded by a family of horospheres for $t>0$. It represents an example of non-uniqueness of nonnegative solutions for the Cauchy Problem.
\end{thm}
We can check in Section \ref{sect.tw} that the special solution, formula  \eqref{gtw.mn}, is not globally bounded nor integrable, as expected for non-uniqueness; actually, its growth for any positive time as the space location goes to infinity is equivalent to $O(s^{1/(m-1)})$, where $s$ is geodesic distance to a given point in $\Hn$.

This is followed in Section \ref{sect.aaas} by the construction of a real asymptotic approximation of fundamental solutions, based on the approximation: $\coth(r)\sim 1$ for $r$ away from 0 in the radial version of the HPME. In this way we get the explicit supersolution (\ref{press.formula}) which will turn out to be a very good approximation of the fundamental solutions for large times. We will also get upper and lower estimates valid for all nonnegative solutions with bounded and compactly supported data.

In the rest of the paper we deal with the rigorous justification of the sharper form of the results. One possible way would be to use the information provided by the supersolutions plus a suitable scaling technique to obtain a family of rescaled problems that converge to the desired profile with the correct rates. This would parallel the direct method that is used in the Euclidean case, \cite{Vaz2003}.

We will use instead the formal equivalence of our problem with an already solved problem to transfer the results from there to the original issue. Thus, in Section \ref{sect.change} we introduce the change of variables that allows to reformulate the evolution in hyperbolic  space in radial coordinates as a Euclidean evolution for a porous medium equation with a weight in Euclidean space (so-called theory of inhomogeneous media).  We then introduce the known results for such a theory, Section \ref{sect.wpme}, and derive the new ones that we need. We combine these two sections in  Section \ref{sect.app} to derive the asymptotic behaviour of the fundamental solution. Not only the solution is estimated in a sharp way, but also the location of the free boundary is.

An important consequence  of the transfer method is that we can describe at the same time the long-time behaviour of the fundamental solutions and also the behaviour of the wide class of weak solutions with  nonnegative integrable, in accordance with the expected role of fundamental solutions in diffusive theories. This is carefully explained in Section \ref{sect.app}, where very detailed versions of the different items of the main theorem \ref{thm.main} are stated and proved for general radial solutions. Less sharp results follow too for data without the property of radial symmetry, but the sharp behaviour in that case is still an open topic.

A drawback of our transfer method is that it does not work in dimension $n=2$. This case is examined in Section \ref{sect.2d}. Consequently, our results are not so sharp.  Note that $n=1$ presents no problem since there is no difference between hyperbolic and Euclidean space then.

The paper closes with some additional information. There are a number of possible extensions of our methods; in Section \ref{sect.pl} we give the basic details about the extension to the $p$-Laplacian evolution equation. Section \ref{sect.comm} contains general comments, extensions and open problems. We conclude by an Appendix that contains some calculations that may be useful.

\medskip

\noindent {\bf Comments on the main results}

\noindent Formulas (\ref{fundsol.hyper.1}) and (\ref{press.formula}) show that the solution takes on a conical shape as a function of $r$ for large times. In  Section \ref{sect.aaas} we show a way of seeing such asymptotic behaviour in hyperbolic space as an outgoing radial traveling wave for a transformed problem, and the resulting formula copies the well-known traveling wave front of the Euclidean PME in a very different context. Let us remark that conical shapes are not very common as profiles of diffusive processes, which are mostly rounded, but they are found sometimes. Let us mention two examples of conical shapes in our experience: (i) the resolution of the avalanche phenomenon in complete blow-up  for reaction diffusion equations of the form $\partial_t u=\Delta u+u^p$, done in a paper with Quir\'os and Rossi \cite{QR2004}; and (ii) the profiles the convection diffusion  model $\partial_t u=\Delta_p u+ |\nabla u|^q$  in the case of critical exponents, studied with Iagar and Lauren\c cot in \cite{ILV2011}, where logarithmic scales are also present. Our approximate analysis below shows that the convective effects are very important, but the diffusive term enters to determine the parameter $\gamma$ which controls the precise penetration of the front.

In view of  formula (\ref{fundsol.hyper.1}) together with the fact that for small times the hyperbolic solutions are a small perturbation of the Euclidean ones given by (\ref{Bar.euc.1})-(\ref{Bar.euc.2}), we conclude that the behaviour of the fundamental solution is very different for small and large times. In particular, it cannot be self-similar. We can actually view the fundamental solution as a smooth connection in function space from the standard Barenblatt function, to which it tends as $t\to 0$, towards the conical wave (\ref{press.formula}) that is the asymptotic limit as $t\to\infty$. Summing up,  curvature combines with size to mark the gradual influence of the curved geometry on the form of the fundamental solution in hyperbolic space. We will  continue the discussion of the results and their consequences in Section  \ref{sect.comm}.


\section{Laplacian and heat equation on hyperbolic space}
\label{sect.2}

The existence of solutions of the heat equation  and the PME  in the hyperbolic space is our concern. We recall that several models can be used to describe $\Hn$ in an explicit coordinate system. For instance, one may realize $\Hn$ as an embedded hyperboloid in $\re^{n+1}$, endowed with the inherited Minkowski metric. It is also possible to use one of the two Poincar\'e realizations, in the sense that topologically one can identify $\Hn$ with the unit ball in $\ren$ or with the upper half-space, each of which endowed with an appropriate metric with the property that the Riemannian distance from any given point to points approaching the topological boundary tends to $+\infty$. Another possible realization is the Klein model. See \cite{Bear83, Chavel84, Gr2009, Rat06, Th97} for a comprehensive discussion. Because of the structure of the isometry group of $\Hn$ it is often convenient to describe the hyperbolic space as a model manifold as described by Green and Wu, \cite{GW79}. On such a manifold, a pole $O$ is given and the Riemannian metric has the form
$$
ds^2 = dr^2 + f(r)^2d\sigma^2,
$$
for an appropriate function $f$, where $r $ is the geodesic distance from the pole $O$ and
$d\sigma^2$ denotes the canonical metric on the unit sphere $\mathbb{S}^{n-1}$. The hyperbolic space $\Hn$ is obtained  making the precise choice $f(r) = \sinh r$.

The theory can be addressed more easily in the setting of radial functions (in the sense we have mentioned above), and this is sufficient for many purposes, like the properties of fundamental solutions. It is known, see for instance \cite{Dav89} and references therein, that the radial part of the Laplacian (more precisely, the Laplace-Beltrami operator) in hyperbolic space has the explicit expression (i.\,e., acting  on radial functions $ u(r)$)
\begin{equation}\label{lapl.radial}
\Delta_{g,rad}u(r) = u''(r) + (n-1)\coth (r) u'(r) =
\frac1{(\sinh r)^{n-1}}((\sinh r)^{n-1}u'(r))'\,.
\end{equation}
 Recall that in such coordinates the volume element is \ $d\mu = (\sinh r)^{n-1} dr d\omega_{n-1}$, where $d\omega_{n-1}$ is the volume element on the unit sphere ${\mathbb S}^{n-1}$. Actually, The differential expression of the Laplace-Beltrami operator $\Delta_g$ in the hyperbolic space with curvature $\kappa<0$ is
$$
\Delta_g u = (\sinh(r/R))^{1-n}\frac{\partial}{\partial r}
\left((\sinh(r/R))^{n-1}\frac{\partial u}{\partial r}\right)
+\frac{1}{R^2\sinh(r/R)^2} \Delta_{\mathbb{S}^{n-1}}u\,,
$$
where $R^2=-1/\kappa$. Reviewing our analysis below will easily convince the reader that there is no essential influence of the particular value of $R$ so that we take without loss of generality $\kappa=-1, R=1$.

\medskip

\noindent {\bf Fundamental solution of the linear heat equation.}
The explicit heat kernel in hyperbolic space, i.\,e., the solution with initial data $U_0(x)=\delta(x)$, is known: in odd dimensions $n=2m+1$ it has the formula
\begin{equation}
U_n(r,t)=\frac{(-1)^m}{2^m\pi^m}\frac1{(4\pi t)^{1/2}}\left(\frac1{\sinh r}\frac{\partial}{\partial r}\right)^m\,e^{-m^2t-r^2/4t}\,,
\end{equation}
 while for even $n$, $n=2m+2$, we have
$$
U_n(r,t)=\frac{(-1)^m}{2^{m+5/2}\pi^{m+3/2}}t^{-3/2t}e^{-\frac{(2m+1)^2}4 t}\left(\frac1{\sinh r}\frac{\partial}{\partial r}\right)^m \int_r^\infty \frac{se^{-s^2/4t}} {(\cosh s - \cosh r)^{1/2}}\,dr\,.
$$

Of course, for $n=1$ we recover the Euclidean solution. For $n=3$ we have
\begin{equation}\label{decay.heat}
U_3(r,t)=\frac{1}{ (4\pi t)^{3/2}}\,\frac{r}{\sinh r}   \,e^{-t-r^2/4t}
\end{equation}
and for $n=2$ the rather complicated expression
$$
U_2(r,t)=\frac{\sqrt2}{(4\pi t)^{3/2}}e^{-\frac{1}4 t} \int_r^\infty \frac{se^{-s^2/4t}} {(\cosh s - \cosh r)^{1/2}}\,dr\,.
$$
Of course, $r$ is the geodesic distance to $O$. For a reference to these formulas see  \cite{GN1998}.

\medskip

\noindent {\bf Comments.} The reader may want to use these formulae to compare with the results of the porous medium case derived in this paper, which are not so explicit but surprisingly simpler in some sense. For instance, it is interesting to compare, say for $n=3$, the decay  (\ref{supest}) of the HPME, $U(0,t)\sim (\log(t)/t)^{1/2}$, with the exponential decay for the heat equation: \ $U(0,t)\sim e^{-t}t^{-3/2}$, that follows from formula (\ref{decay.heat}). The alternative between power decay for PME versus exponential decay for the linear HE occurs also in Dirichlet problems posed in bounded domains of Euclidean space.

Moreover, it is easy to see that formula (\ref{decay.heat}) for the heat equation predicts an approximate separate variables behaviour at all istances $r=o(t)$, while this is false for the HPME where the term $(\gamma \log t -r+b)_+$ in formula (\ref{fundsol.hyper.1}) breaks this behaviour at distances of the order of $\gamma \log t$.

The description of asymptotic behaviour in terms of special solutions of the separate-variables type or the self-similar type is a recurrent fact in the theory of partial differential equations, \cite{Barbk96, Sachdev, S4, VazSm2006}, and it happens for the standard  examples of linear and nonlinear diffusion in the Euclidean setting. It does not happen, at least in an open way, in the problems we are considering in hyperbolic space, which adds complication to the study and interest to the mathematical analysis.

\section{ The HPME.  Basic ideas and first estimates}
\label{sect.3}

The theory of existence and uniqueness of solutions of the HPME (\ref{hpme}) with standard classes of initial data is not the main concern of this paper, so we will be rather brief. It can be done in the framework of the theory of maximal monotone operators in the space $H^{-1}(\Hn)$, extending the ideas of  Brezis's paper \cite{Brez71} from the Euclidean space to the manifold setting. This is first done in the book \cite{Vbook} (Chapter 11) for manifolds with nonnegative Ricci curvature, and in \cite{BGV08} for fast diffusion in more general manifolds, see also \cite{Dek08, GM2013}.
The changes to cover the PME on the hyperbolic space are minor. We also refer  to \cite{Pu12} where the existence of solutions of a more general problem is done and the basic propagation properties established. In particular, it is clear that the $L^p$ norms of the solutions, $1\le p\le \infty$,  are non-increasing in time, that the solutions form a family of $L^1$-contractions,  that the Maximum Principle works, and that the mass conservation law holds
$$
\int_{Hn} u(x,t)\,d\mu(x)=\int_{Hn} u(x,0)\,d\mu(x).
$$
At this stage we do not need any assumption of radial symmetry, we may just start with data
$u_0\in H^{-1}(\Hn)\cap L^1(\Hn)$. Here we are only interested in nonnegative solutions. We want to concentrate mainly, though not exclusively, on the fundamental solutions. As we said above, there is no loss of generality in considering only the case of initial data of unit mass.

We will show below (Section \ref{sect.change}) that the case of nonnegative and radial data and solutions  can be reduced to the equivalent study of a porous medium flow with a weight posed in the Euclidean space, where the questions of existence, uniqueness and main properties are already solved. We will accept this convenient fact to proceed to the goal without much delay to keep the exposition within a reasonable extension.


\subsection{Properties }

The behaviour for small times does not have any problem if we start with compactly supported data, since the finite speed of propagation implies that the support will be finite and expanding for positive times. In this way, it is easy to pass from continuous radial initial data to the existence of the fundamental solution.

Moreover, it is easy to prove that radially decreasing data produce solutions that are radially decreasing in space for every time, and this is also true in the limit for the fundamental solution.

The uniqueness of the fundamental solutions in the setting of radial solutions follows from the same result for the PME posed in an inhomogeneous medium, according to the equivalence of problems that we will establish in Section \ref{sect.change} and ff. Without that restriction, it can be obtained by following the proof for the PME in the Euclidean case that consists in passage to the so-called dual equation, as done by Pierre in a classical paper \cite{P82} in the Euclidean case,  and recently extended by the author \cite{Vz2014} to the fractional porous medium equation in the case initial data a Dirac measure, and then in \cite{GMPu2013} where the initial data are positive finite measures. This implies some technical work that is not the object of this paper and will appear elsewhere.  Finally, we recall that an alternative proof of uniqueness of fundamental solutions in the case of the PME and the $p$-Laplacian equation in Euclidean space was proposed by Kamin and the author in \cite{KV88}, but it uses heavily the scaling groups and it does not seem well suited to be adapted to the hyperbolic space setting.

We now proceed with preliminary estimates to guess the sizes that we will deal with.

\medskip

\noindent {\sc $L^\infty$ estimates.} (i) In order to obtain a sup bound from above, we only have to compare  with PME, since a simple inspection of the equation shows that radially decreasing solutions are subsolutions of Euclidean PME
$$
U(x,t)\le Ct^{-\alpha_e}, \qquad \alpha_e=\frac1{m-1+(2/n)}
$$
This estimate should be very good for small times, but it is not good later when the curvature of the space begins to have a strong influence, as we will show.

\noindent (ii) To get a bound from below, we  compare with the Dirichlet Problem posed in a ball of  radius $R$ and we get
$$
U(r,t)\ge F(r)t^{-1/(m-1)}\,,
$$
for $0<r<R$ and $t>0$, which  means that the fundamental solution decreases no more than the  power-rate that we are familiar with as the worst case in the Euclidean theory (i.\,e., the bound that holds for all dimensions). The proof uses the method of subsolutions and we leave the easy detail to the reader.

\medskip

\noindent {\sc Propagation.} (i) It is very simple to prove that there is eventual penetration of the support into all the space. Otherwise, there would be  localization in a ball for all time, which would make the problem look more or less like a porous medium in a bounded Euclidean domain, and this is impossible. This is only a rough idea, but enough in view of the detailed information that follows below.

 \noindent (ii) More quantitatively, if the solution is supported in the ball of radius $R(t)$, then a bound from above for the $L^\infty$ norm of the solution of the form $\|u(t(\|_\infty \le Ct^{-\alpha}$ and conservation of mass  imply that
$$
\sinh R(t)\ge C(t^{\alpha/(n-1)})\,,
$$
which means that for all large $t$ we have
$$
R(t)\ge  \frac{\alpha}{n-1}\log t +c_1\,.
$$
This is the first indication of the actual logarithmic growth that will be rigorously proved in the paper. Once the correct $\alpha=1/(m-1)$ is found, the last estimate will turn out to be correct too.

  \noindent (iii) On the other hand, the same comparison with the Euclidean case shows that the free boundary has the upper  bound $R(t)\le C t^{\beta}$, $\beta=1/(n(m-1)+2)\,.$
It will turn out that this estimate is very unrealistic.

\medskip

\noindent (iv) Bulk estimate. A good upper estimate for $R(t)$ is not so easy. In order to show an upper estimate of logarithmic growth we consider $R(t)$ the distance of half-height values for $U(\cdot,t)$. Then the mass there is of order \ $c\,t^{-1/(m-1)}V(R)$,
and since mass is conserved and  $V(R)\sim (\sinh R)^{n-1}$ we get
$$
\sinh R(t)\le Ct^{1/(m-1)(n-1)}
$$
This can be done for every positive ratio, not only $1/2$. The limit of ratios going to zero will be the free boundary estimate, which needs another analysis.

\medskip

\noindent (v) Tail analysis. This applies to all radially symmetric decreasing solutions in $L^1(\Hn)$ and comes from conservation of mass. It reads as follows
$$
U(r)\le \|u_0\|_1/V(r)= c\|u_0\|_1/(\sinh r)^{n-1}\sim c\|u_0\|_1 e^{-(n-1)r}
$$
(the last formula for all large $r$) which is a first estimate, uniform for all times. It is not very accurate when the profiles are actually compactly supported and bounded.

\medskip

\noindent $\bullet$ {\sc Universal $\partial_t u$ estimate.} There is an a priori estimates based on homogeneity, originally due to  B\'enilan et al. \cite{AB79, BC81} in the Euclidean case and valid for all nonnegative solutions.  We have
\begin{equation}\label{beni.est}
\partial_t u\ge -\frac1{(m-1)t}\,.
\end{equation}
It is now a standard argument, see \cite{Vbook}, that this implies that any point where the solution is positive for some $t>0$ will be positive for $t'>t$, a property usually called retention, that implies that the support of the solution is monotonically non-increasing family of sets in $\Hn$.

\medskip

\subsection{Short time limit of the fundamental solutions}

In order to analyze the behaviour for small times we take a fundamental solution $U(r,t)=U(r,t;M)$  of the HPME  with hyperbolic mass $M>0$ and perform the typical Euclidean PME rescaling
$$
U_k(r,t)=k^{n\beta}U(k^{\beta}r,kt)
$$
with scaling parameter $k>0$ (that will be small). Then $U_k$ satisfies the equation
$$
\partial_t U_k=\partial^2_{rr}(U_k^m)+(n-1)k^{\beta}\coth(k^\beta r)\,\partial_{r}(U_k^m)
$$
Since $\coth(r)\ge 1/r$ solutions of this equation with $\partial_{r} U\le 0$ are subsolutions of the Euclidean PME. Moreover,
$$
M=c_n\int_0^\infty U(r,t)\,\sinh^{n-1}(r)dr\ge c_n\int_0^\infty U(r,t)\,r^{n-1}dr\,,
$$
and since the support of $U(r,t)$ shrinks to the origin as $t\to 0$ we have
$$
c_n\int_0^\infty U(r,t)\,r^{n-1}dr \to M \quad \mbox{as \ } \ t\to 0.
$$
We conclude that $U(r,t)$ lies below the Euclidean fundamental solution with the same mass, described in (\ref{Bar.euc.1})-(\ref{Bar.euc.2}), and the same happens for all $U_k$.

Next using the fact that the function $f(r)=r\coth(r)$ is increasing in $r$, we conclude that
the sequence of functions $U_k(r,t)$ is monotone decreasing as $k$ decreases. Due to the a priori bound by the Euclidean solution, we can pass to the limit in the weak formulation as $k\to 0$. Since we have
$k^{\beta}\coth(k^\beta r)\to 1/r$ as $k\to 0$ for positive $r$ we get as limit of the $U_k$ a solution of the Euclidean PME for $r>0$. Since the support goes to the origin as $r\to 0$ and the (Euclidean) mass is $M$, by uniqueness we know that it is the Barenblatt solution, that we call $B(r,t;M)$ for definiteness. Therefore, setting $t=1$ in the limit we get
$$
\lim_{k\to 0}k^{n\beta}U(k^{\beta}r,k)=B(r,1;M)
$$
Changing variables ($k$ into $t$ and $K^{\beta}r$ into $s$) and using the scale invariance of $B$ we get the following result

\begin{thm} If $U(r,t;M)$ is the fundamental solution of the HPME and $B(r,t;M)$ is the Barenblatt solution for the Euclidean PME, both with mass $M>0$ w.r.t the respective volume measure, then
\begin{equation}
\lim_{t\to 0}t^{n\beta}|B(r\,t^{\beta},t;M)-U(r\,t^{\beta},t;M)|=0
\end{equation}
with pointwise limit for all $r>0$. Moreover, $|B(\cdot,t;M)-U(\cdot,t;M)| \to 0 $ in $L^1(\ren)$.
\end{thm}


\section{An explicit typical example of finite propagation}\label{sect.tw}

The construction of explicit solutions plays an enormous role in nonlinear theories. Sometimes the explicit solutions that we can construct do not belong to the classes of solutions that we are interested in for physical or mathematical reasons, but nevertheless they usually offer valuable intuitions of qualitative or quantitative aspects of the theory, since much of the mathematical analysis is based on local behaviour.

In the case of the PME in Euclidean space this occurs with the {\sl traveling wave solutions} given by the formula $p(x,t)=c(ct-x_1)_+$, $c>0 $ ($p$ is the pressure). They exhibit many of the properties of the actual propagation front of the practical solutions, for instance solutions with compactly supported and integrable data. Even if the global properties are different from the class of integrable solutions, they are good models at the local level. So we may wonder what could a similar construction in hyperbolic space. We propose here a solution to that question.

\noindent $\bullet$  Let us examine first the simplest case,  $m=2$ and $n=2$, for clarity.  The preferred representation of the hyperbolic space will be  Poincaré's upper half-plane representation where $\mathbb{H}^2$ is identified with $\{(x,y): x\in\re, y>0\}\subset \re^2$, the metric is given by
$$
ds^2=y^{-2}(dx^2+dy^2)\,,
$$
and the Laplacian operator has  the expression
$$
\Delta_{\mathbb{H}^2} u=y^2 \left(\frac{\partial^2 u}{\partial x^2} +  \frac{\partial^2 u}{\partial y^2}\right)\,.
$$
It is then easy to see that for every $c>0$ the function
\begin{equation}\label{gtw.22}
U(x,y,t)=\frac{(\log(cty))_+}{2t}
\end{equation}
is a weak solution of the HPME, in the class of self-similar solutions, since it has the form $U({\bf x},t)=t^{-\alpha_1}F({\bf x}t^{-\alpha_2})$, where ${\bf x}=(x,y)$. Here $\alpha=1$ and $\alpha_2=-1$. The equation for $F(\xi)$ with $\xi= yt$ is then
$$
-F+\xi F'=\xi^2 (F^2)''\,.
$$
Indeed, a simple calculation shows that the equation is satisfied everywhere in the positivity set $\{(x,y,t): U>0\}$. On the other hand, in order to obtain a global weak solution we have to check the conditions at the free boundary,  which is the set
\begin{equation}\label{fb.gtw}
\Gamma=\{(x,y,t): y=\frac{1}{ct}\}
\end{equation}
and is planar front in the Euclidean graph of this space representation. It is well known in the PME theory that a condition like the continuity of $U$ plus the space differentiability of $U^m$ at the free boundary (FB for short) implies that the weak formulation is satisfied, so we will not repeat the easy calculation. Moreover, we can obtain a sharper condition which is Darcy's law.  Let us review that important point: using (\ref{fb.gtw}) the normal speed of propagation of the  free boundary is
$$
\frac{ds(t)}{dt}=\left|\frac{ds}{dy}\,\frac{dy(t)}{dt}\right|=\frac{1}{y(t)}\,\frac1{ct^2}=\frac1{t}\,,
$$
where $s(t)$ measures geodesic distance from some given origin located at $y=1$. We immediately see that Darcy's law holds in the sense that  at the FB the geometrical speed $ds(t)/dt$ equals the speed understood in the standard internal way, i.\,e., the physical derivation says that the speed vector equals minus the gradient of the pressure $P=2U$ (cf. \cite{Vbook}, Section 2.1) which amounts to
$$
-\partial_s P=-2\partial_s U=-2\partial_y U/(ds/dy)=1/t.
$$
In fact, this internal speed is constant everywhere at a given time, which justifies the name of (generalized) traveling wave. Standard TWs have the form $F(x_i-ct)$ with constant $c$, and this is not the case here, hence the word `generalized'. See further comments in the final section.

The example is a close  equivalent to the usual traveling waves of the PME which have however the standard traveling wave form. Like them, our TW has a plane front lines. Here these lines,  $y=$ constant, are part of the family of the special curves called {\sl horocycles} in hyperbolic plane. In the representation as Poincaré's disk they are a family of (Euclidean) circles tangent to the unit circle $\|x\|=1$ at any horizon point, say $(-1,0)$, so that the solution can be  seen as a wave advancing with the farthest tip on the $x_1$ axis and tending to $(1,0)$ with rate $x_1(t)= ct/(ct+1)$. Using the formulas for the geodesic distance from the origin we get $s(t)\sim \log(t)+b$, where $b=\log c$. This is also in perfect agreement with point (iv) of our Theorem \ref{thm.main}.

\medskip

\noindent$\bullet$ {\bf General parameters.} The calculation for general $m>1$, $n\ge 2$ goes along the same lines with suitable changes in the exponents that bear in mind that the pressure function, given by $P=(m//m-1)U^{m-1}$, should look similar in all cases. The explicit self-similar function is now given by
\begin{equation}\label{gtw.mn}
U(x,y,t)^{m-1}=a\frac{(\log(ct^{\gamma}y))_+}{t}
\end{equation}
with $a=1/m(n-1)$ and $\gamma=1/(m-1)(n-1)$, which coincide with the values found in the main theorem. We leave to the Appendix the verification of the equation in the positivity set; we also have to check the Darcy condition on the free boundary, which is done as before. This one is given by the equation $y(t)=1/(ct^{\gamma})$, and the normal speed of the free boundary is $ds/dt=\gamma/t$, with leads to the expression $s(t)\approx \gamma \log(t)$. For the reader's convenience we recall that the Laplacian operator in this representation of $\Hn$ is given by the formula
$$
\Delta_{\mathbb{H}^n} u=y^2 \left(\frac{\partial^2 u}{\partial x^2} +  \frac{\partial^2 u}{\partial y^2}\right)-(n-2)\,y\,\frac{\partial u}{\partial y}\,,
$$
where we use the analyst's convention on the sign of the Laplace operator.

\medskip

\noindent $\bullet$ {\bf Hidden traveling wave.} In order to understand the structure of the explicit solution in more classical terms it is convenient  to use a renormalization that is common in Dirichlet problems in bounded domains of Euclidean space:
\begin{equation}\label{renorm1}
w(x,\tau)=t^{1/(m-1)}u(x,t), \qquad \tau=\log(t)
\end{equation}
(note the new logarithmic time scale). This transforms the HPME into the reaction-diffusion equation
\begin{equation}\label{renorm2}
w_\tau=\Delta_g (w^m)+ cw, \qquad c=1/(m-1)\,.
\end{equation}
 In that setting the explicit solution $ U$ of formula (\ref{gtw.mn}) can be written as
\begin{equation}
P(x,y,t)=\frac{m}{m-1}{W}^{m-1}(x,y,t)=\gamma(\gamma\tau+\log(yc))_+\,,
\end{equation}
which is an outgoing radial traveling wave when we take variable $r=-\log(y)$, and then the wave speed is just $\gamma$ everywhere that ${W}>0$. Actually, the  reader will not fail to recognize the well-known formula for the  traveling wave of the PME in Euclidean space in pressure variable, a lucky formal coincidence!\footnote{But notice that in the Euclidean  traveling wave $\gamma$ is not fixed, it can have any positive value.} Notice also that ${P}^{m-1}(r,t)$ grows linearly in time everywhere.
Note finally the growth rate $U \sim \log(y)^{1/(m-1}$ as $y\to\infty$ for fixed $t>0$.

\medskip

\noindent $\bullet$ {\bf Initial trace and non-uniqueness.} Let us now address the contents of Theorem \ref{thm1.2}. Indeed, the constructed solution has a very definite support that forms a certain neighbourhood of the particular boundary point at infinity represented by the direction $y\to +\infty$ in our half-space representation of hyperbolic space, and this support shrinks as time $t\to0$, so that for every compact subset of $\Hn$ there is a time $t(K)$ such that $U(x,t)=0$ if $x\in K$ and $0<t<t(K)$. We conclude that the initial trace is zero in any (current) sense of local convergence in $\Hn$.

Since solution \eqref{gtw.mn} is non trivial and the trivial function is also a solution with zero initial trace, this creates a non-uniqueness phenomenon in the form of ``a wave coming from infinity'' (more precisely, from a point in the infinite boundary of hyperbolic space in the Poincar\'e disk/ball representation). Such phenomenon never happens for nonnegative solutions of the heat equation or the PME in the Euclidean space $\ren$ where it is known that non-uniqueness may happen only for highly oscillating solutions. In the case of the HE the non-trivial oscillating solution with zero trace is due to Tychonoff \cite{Ty35}, 1935. Self-similar oscillating solutions with zero initial trace for the PME were constructed much later by Vázquez \cite{Vaz90}.

\medskip

\noindent {\bf Further Comments.} (1) Formula \eqref{gtw.22} produces a unique special solution of that form for $n=2$ and $m=2$. Passing from the half-plane representation to the Poincar\'e disk representation we immediately see that the formula  not only leads to one solution but to a family of solutions, more precisely, one solution starting from every horizon point. The same happens for $n\ge 3$ and all $m>1$.

The problem is posed whether there are more special solutions of a similar kind and what they look like. For instance we can start special solutions from two horizon points $A_1$ and $A_2$: as long as they have disjoint supports the sum of the two is still a solution of the HPME.

(2) By considering  in this representation for $n=2$ more general solutions that still are $x$-independent,  we are faced with  the 1-D equation: $y^{-2 }  u_t=(u^m)_{yy}$, which leads us to the theory of nonlinear diffusion in inhomogeneous media with inverse square weights, a remarkable limit case in diffusion theory. This idea cannot pursued further here, but it will lead to very interesting consequences in Sections \ref{sect.change}--\ref{sect.app} and in dimensions $n>2$.

\section{Asymptotic analysis via approximate solutions}\label{sect.aaas}

The next step of our analysis is to get precise asymptotic estimates for the fundamental solutions and other radial solutions with integrable data in the form of actual profiles and rates. To start with, in order to estimate the actual rate of divergence of the free boundary location, we may use the intuitive idea that the hyperbolic metric implies an ample space in transversal directions  as compared with the Euclidean one, and this will lead to a asymptotic estimates like a Euclidean PME with dimension $n\to  \infty$, so that we expect the rate of increase to be less than any power of time. The example of the previous section indicates that we will end up with some logarithmic expression of time.

\subsection{Explicit solution of the approximate equation}

In order to make this intuitive idea precise in a quantitative way we present next an analysis based on a sharp approximation that works as follows: since the solution will be spread with time over a large part of the space we may replace for all practical purposes  the coefficient $\coth(r)$ by 1 in the radial expression for the Laplacian  (\ref{lapl.radial}), and consider exact solutions of the following approximate equation in radial coordinates
\begin{equation}\label{eq.apprx}
\partial_t u=(u^m)_{rr}+(n-1)(u^m)_r\,.
\end{equation}
After the preliminary analysis and some other evidence, we propose to try a solution of this  equation of the form
\begin{equation}
{\widetilde U}(r,t)= t^{-1/(m-1)}F(\gamma \log t-r+b)\,.
\end{equation}
If $\xi=\gamma \log(t)+b-r$, then $F(\xi)$ must then satisfy
\begin{equation}\label{profileF}
-\frac1{m-1}F(\xi)+\gamma F'(\xi)= (F^m)''+(n-1)(F^m)'\,,
\end{equation}
where primes denote derivatives with respect to $\xi\in\re$. In this way the time dependence of equation (\ref{eq.apprx}) is eliminated (also the constant $b$ disappears). We  now propose for the profile function $F$ the choice \ $F(\xi)=(a\xi)_ +^{1/(m-1)}$, which copies the example of Section \ref{sect.tw} for a different kind of coordinates. The coefficients $a$ and $\gamma$ have to be determined. Let us make the calculations with some detail  this time. Inserting the ansatz into   (\ref{profileF}) we get the conditions to be satisfied, at least for $\xi \in (0,\gamma \log(t)+b)$:
\begin{eqnarray*}
\displaystyle &-\dfrac{a^{1/(m-1)}}{m-1}\xi^{1/(m-1)} + \dfrac{\gamma a^{1/(m-1)}}{m-1}\xi^{1/(m-1)-1}\\[10pt]
=
&-\dfrac{(n-1)ma^{m/(m-1)}}{m-1}\xi^{1/(m-1)} + \dfrac{ma^{m/(m-1)}}{(m-1)^2}\xi^{1/(m-1)-1}\,,
\end{eqnarray*}
and after simplifying
$$
-\xi +\gamma = -(n-1)ma \xi+\frac{ma}{m-1}\,.
$$
This is satisfied for the precise values
$$
a=  \frac{1}{m(n-1)}\,, \quad \gamma=\frac{ma}{m-1}=\frac{1}{(m-1)(n-1)}\,.
$$
We see that  $b\in \re$ is a free parameter (corresponding to space translation invariance). This constant can be used to adjust the mass of the solution at any given time. In order to complete the analysis of the approximate solution we point out that the equation is trivially satisfied in the open set where $r> \gamma \log(t)+b$ since there $F=0$ (so-called empty set). Finally,  we need to be aware that there is a problem of regularity of $\widetilde U$ at the free boundary  $\{(r,t): \ r= \gamma \log(t)+b\}$, but we may argue as in the standard PME, by pointing out  that ${\widetilde U}^m$ is a $C^1$ function for all $r>0$, and this is enough to imply the weak solution condition there.

 This approximate solution
 \begin{equation}\label{approx.sol}
 {\widetilde U}(r,t)=\frac{(a(\gamma \log(t)+b-r))_ +^{1/(m-1)}}{t^{1/(m-1)}}
 \end{equation}
 agrees with the main result of Theorem \ref{thm.main} and allows to calculate $a$ and $ \gamma$ in a very clear and easy way.
We will prove in Section \ref{sect.app} that ${\widetilde U}(r,t)$  is a very good approximation of the actual solutions as $t\to\infty$ once $b$ is conveniently chosen. We ask the reader to check for himself that for small $t$ this approximate solution does not resemble at all the  actual solutions of the PME in hyperbolic space, which look like Euclidean.

Note finally that we can construct a whole family of approximations by time displacement, ${\widetilde U}(r,t;\tau)={\widetilde U}(r,t+\tau)$, $\tau\in\re$.

\medskip

\noindent $\bullet$ {\bf Supersolutions.} Since the approximate solution ${\widetilde U}$  satisfies ${\widetilde U}_r\le 0$ and  $\coth (r)\ge 1$, it is immediate that for $r>0$ we have
$$
\partial_t {\widetilde U}\ge \Delta_g {\widetilde U}^m\,,
$$
which means that it ${\widetilde U}$ is an actual supersolution of the original equation $\partial_t U= \Delta_g U^m$  on $\Hn-\{0\}$. We ask the reader to check that it is also a supersolution at  $r=0$ in the weak sense. Since we have already justified the situation at the free boundary, we are allowed to use the  standard comparison results between solutions and supersolutions for parabolic equations.

We immediately conclude from such comparison that the estimates from above contained in Theorem \ref{thm.main} hold for all solutions with compactly supported initial data, both as regards the size of the solution and also the support, but of course we expect the constant $b$ not to be exact, but approximate from above. In particular, it is true for any solution  with compact support that the free boundary lies in a ball with radius
\begin{equation}\label{sup.fb}
R_+(t)\le \gamma (\log t)+b
\end{equation}
with $\gamma=1/(m-1)(n-1)$ and some $b$ that depends on the solution. It is also proved that
\begin{equation}\label{sup.li}
\|u(\cdot,t)\|_\infty \le t^{-1/(m-1)}(c_1\log(t)+c_0)^{1/(m-1)}\,, \quad c_1=a\gamma.
\end{equation}
In these estimates we are not assuming radial symmetry.

\medskip

\noindent $\bullet$  {\bf Free Boundary condition.} It is a relation that links the speed of propagation of the free boundary $R'(t)$ with the internal velocity of propagation associated with the porous medium equation, ${\bf v}=- (m/(m-1)\nabla u^{m-1}$. It holds for our approximate solution in the form
\begin{equation}
R'(t)=-\frac{m}{m-1}\,\left.\partial_r ({\widetilde U}^{m-1})\right|_{r=R(r)-}=\frac{\gamma}{t}.
\end{equation}
Actually, ${\bf v}$ is minus the gradient of the pressure discussed in the Introduction, a fact that reflects the underlying physics (Darcy's law) and explains its importance in studies of support propagation.

\medskip

\noindent $\bullet$ {\bf Hidden traveling wave.} As in Section \ref{sect.tw} we can improve the understanding of the structure of the asymptotic approximation using the renormalization
\begin{equation}\label{renorm1b}
w(x,\tau)=t^{1/(m-1)}u(x,t), \qquad \tau=\log(t)
\end{equation}
which leads to the reaction-diffusion equation $w_\tau=\Delta_g w^m+ cw$ with $c=1/(m-1).$
Then the asymptotic approximate solution $\widetilde U$  can be written as
\begin{equation}
\frac{m}{m-1}{\widetilde W}^{m-1}(r,t)=\gamma(\gamma\tau-r+b)_+\,,
\end{equation}
which is an outgoing radial traveling wave, whose speed is just $\gamma$ everywhere that ${\widetilde W}>0$. Comparison with the standard traveling wave of the PME in Euclidean space, written in pressure variable, is also natural. Moreover, ${\widetilde W}^{m-1}(r,t)$ grows linearly in in time everywhere, and it has an exact conical spatial shape.

\medskip

\noindent $\bullet$  {\bf Mass calculation.}
We may calculate the hyperbolic mass of the approximate solution
$$
M(t)=\int_{\Hn} U(x,t)d\mu(x)=c_nt^{-1/(m-1)}\int_0^\infty a(\gamma  \log t-r+b)_+^{1/(m-1)}(\sinh r)^{n-1} dr\,.
$$
Using the change $r=\gamma\log(t)+b-\xi$ and the approximation $\sinh(r)\sim (1/2)e^t$
and noting that $\gamma (n-1)=1/(m-1)$ we get
$$
{\widetilde M}(t)\sim d e^{b(n-1)}\int_0^{\gamma  \log(t)} \xi^{1/(m-1)}e^{-(n-1)\xi}d\xi
$$
and for $t$ very large
$$
{\widetilde  M}(t)
\sim d e^{b(n-1)}\int_0^{\infty} \xi^{1/(m-1)}e^{-\xi(n-1)}d\xi=K(n,m)e^{b(n-1)}\,.
$$
where $d=d(m,n)$. Since ${\widetilde U}$ is a supersolution and the mass in conserved for solutions, then ${\widetilde  M}(t)$ will be (slightly) increasing with time (this could be checked by direct calculation). However, the limit is precisely the asymptotic mass that comes out for the fundamental solution in the main theorem, as we will see.

\subsection{\bf Construction of a subsolution}

An upper bound for  $\coth(r)$ is $\coth(r)\le (1+r)/r$ for all $r>0$.
We then have $\coth(r)\le k$ for $r\ge 1$ and $\coth(r)\le k/r$ for $0<r<1$ for some  $k >1$. Using this information, we will try to  get a subsolution of our problem by working again with approximate equations, using now inequalities in the other direction. The analysis is not so simple, we have to work separately in the two regions and use matching to join them.

\medskip

$\bullet$ For $r\ge 1$ we propose a subsolution of the  same form as before (since it will be proved that this form is almost exact), but now we write
\begin{equation}
{\widetilde U}(r,t)= t^{-1/(m-1)}F(\gamma \log t-r), \qquad F=(a\xi+b)^{1/(m-1)},
\end{equation}
with $\xi=\gamma\log t-r$, and coefficients $\gamma$, $a$ and $b$ have to be determined again, they will be a bit different from the previous ones. The  subsolution condition reads
$$
-\frac1{m-1}(a\xi+b)^{1/(m-1)} + \frac{\gamma a}{m-1}(a\xi+b)^{1/(m-1)-1}\le
$$
$$
-k\frac{(n-1)ma}{m-1}(a\xi+b)^{1/(m-1)}
+ \frac{ma^2}{(m-1)^2}(a\xi+b)^{1/(m-1)-1}\,.
$$
After simplifying
$$
-(a\xi+b) +\gamma a\le  -k(n-1)ma(a\xi+b)+\frac{ma^2}{m-1}.
$$
which is satisfied if
$$
1 \ge k(n-1)ma, \quad \gamma a-b\le -k(n-1)mab+\frac{ma^2}{m-1}.
$$
The system is satisfied for instance by the values
$$
a=  \frac{1}{km(n-1)}\,,\qquad \gamma\le \frac{ma}{m-1}\,,
$$
and $b$ is then a free parameter.  But we can also use $a<a_0=  1/(km(n-1))$ and then the condition on
$\gamma$ does not change if $b\ge 0$. Observe that this value of $a$ can be obtained by using
the previous calculations for an exact approximate solution, but posed in a higher space dimension, $n'>n$ (so that $n'-1\ge k(n-1)$).

\medskip

$\bullet$ We have to complete this ``outer analysis'' with a subsolution in the ``inner region'' $0<r<1$. Compared with the outer region, which has size of the order of $\log(t)$ in space, this is a ``thin region''. Both have  to be matched to order 1 of regularity to get a global subsolution function.

We know that $\coth(r)\le k/r$, $k >1$, for $r\le 1$. We take a function of the form
$$
U^m=C_1(t)-C_2(t)r^2
$$
as a candidate subsolution in this region. The  subsolution condition inside the region takes the form of the differential inequality
$$
\frac1m(C_1(t)-C_2(t)r^2)^{(1-m)/m}(C_1'-r^2C_2')\le 2C_2(t)(-1-(n-1)(1+r))\le -2K_1C_2(t).
$$
where we have used $\coth(r)\le (1+r)/r$.  and  $K_1$ denotes a constant  that depends on $n$ and $m$, while $t$ is large and $r\in (0,1)$. We can take $K_1=2n-1$. Besides, we need matching conditions  between the inner and outer formulas at $r=1$ for $U^m$ and $\partial_r U^{m}$:
$$
C_1(t)-C_2(t)=A^m, \quad -2C_2(t)=-mA^{m-1}B\,,
$$
where
$$
A:=t^{-1/(m-1)}(a(\gamma\log t-1)+b)^{1/(m-1)}\approx t^{-1/(m-1)}(a\gamma\log t)^{1/(m-1)}
$$
$$
B=\frac{a}{m-1}t^{-1/(m-1)}(a(\gamma\log t-1)+b)^{1/(m-1)-1}\approx \frac{A}{(m-1)\gamma\log t}
$$
with approximations for  large $t$. We also have
$$
C_2(t)\approx A^{m-1}B\sim \frac{A^m}{(m-1)\gamma\log t}, \qquad C_1(t)\sim A^m(t)\,.
$$
Let us now check whether the subsolution inequality is satisfied. It  reads
$$
(C_1'-r^2C_2')\le -2mK_1 C_2\,(C_1-C_2r^2)^{(m-1)/m}.
$$
Since $C_i'(t)\approx -C_i(t)/(m-1)t$, $i=1,2$,  a sufficient condition is
$$
C_1\ge 2m(m-1)K_1 t C_2(t)C_1(t)^{(m-1)/m}.
$$
so that $C_1(t)^{1/m}\ge 2m(m-1)K_1 t \, C_2(t)$, i.\,e., $A(t)\ge 2mK_1 {A(t)^m}/\gamma{\log t}$,
hence
$$
A(t)^{m-1}\le K\gamma\log(t)/ t, \quad K=\frac{1}{2mK_1}.
$$
 We expect this to be  true with a suitable choice of the coefficients. Actually, in view of the value of $A$ obtained above, we only need
$$
a\le K(m,n)\,,
$$
that can be obtained in view of analysis at the end of the previous item (by using the exact approximate formula in a higher dimension). This implies also a bound for the admissible value of $\gamma$.

\medskip

\noindent $\bullet$ {\bf Consequences.} Using these results and the Maximum Principle, we conclude that the solutions of the HPME with bounded, nonnegative and compactly supported initial data have a size estimate of the form
\begin{equation}\label{subsupest}
\|u(\cdot,t)\|_\infty^{m-1}\ge c_1\frac{\log (t)}{t}\,,
\end{equation}
and also the minimal radius of the free boundary will be
\begin{equation}\label{sub.fb}
R_-(t)\ge c_2\log(t)\,.
\end{equation}
for some constants $c_1$, $c_2$. These lower values complement the upper estimates (\ref{sup.fb}) and (\ref{sup.li}). In all these estimates we need not assume radial symmetry.


\section{Change into Euclidean diffusion with weight}
\label{sect.change}

The proof of the sharp results on asymptotic behaviour will be implemented in the next sections in dimensions $n\ge 3$  and for radial solutions by means of a very interesting transformation into a weighted equation in radial variables with underlying Euclidean metric, and the resulting weight looks like $1/r^2$ as $r\to \infty $,  which is a critical power for that theory. We will then use the precise asymptotic analysis for the PME acting on an inhomogeneous Euclidean medium given by a certain density function $\rho$. The behavior of $\rho$ plays a key role in such an analysis.

Let us proceed. We have the equation for $u=u(r,t)$ in hyperbolic space $\Hn$
\begin{equation}\label{eq.hpme}
\partial_t u=\Delta_{rad}(u^m )\,, \qquad \mbox{where} \
\Delta_{rad}\,u(r) =
\frac1{(\sinh r)^{N-1}}((\sinh r)^{N-1}u'(r))'\,.
\end{equation}
and $0<r<\infty$. We take dimension $n\ge 3$. We  want to transform the equation by a change of variable $s=s(r)$ into an equation for ${\widehat u}(s,t)= u(r,t)$ of the form
\begin{equation}\label{eq.wpme}
    \rho(s)\,\partial_t {\widehat u}=\Delta_s {\widehat u}^m\,,
\end{equation}
where $\Delta_s$ is the Euclidean Laplacian in $n$ dimensions (and radial coordinate $s=|x|$).  Using the differential expressions, we easily see that the  change rule
\begin{equation}
\frac{ds}{s^{n-1}}=\frac{dr}{(\sinh r)^{n-1}}
\end{equation}
performs the transformation. Integrating, we see  that for $r\sim 0$ we have $s(r)\approx r$, while for $r\gg 1$ we get $s(r)\sim (\sinh r)^{(n-1)/(n-2)}$ and more precisely
\begin{equation}\label{rtos}
s(r) \approx c(n)\,e^{(n-1)r/(n-2)}
\end{equation}
The constant is $c(n)=((n-1)/(n-2))^{1/(n-2)}2^{-{(n-1)/(n-2)}}$. In particular for $n=3$ we get
$s(r)\sim \frac12 e^{2r}$.

After some easy computations we pass  from the original equation (\ref{eq.hpme}) to the Euclidean equation (\ref{eq.wpme}) with weight:
\begin{equation}\label{rho.def}
\rho(s)=\frac{(\sinh r)^{2(n-1)}}{s^{2(n-1)}}\,.
\end{equation}
We conclude that $\rho(s)\approx 1$ for small $s$, while $\rho(s)\approx c_1(n)s^{-2}$ for large $s$.
Note that
$$
(\sinh r)^{n-1}dr= \rho(s)s^{n-1}ds\,,
$$
as needed by the  conservation of mass laws of both equations. It is important to realize that we get a weight that is not a pure power function.

\medskip

\noindent {\bf Comment.} If we perform the same transformation starting from the approximate solution
$$
\partial_t u=(u^m)_{rr}+(n-1)(u^m)_r=e^{(1-n)t}(e^{(n-1)t}(u^m)_r)_r\,,
$$
we get the transformation $ds/s^{n-1}=dr/e^{(n-1)t}$, so that in this case
$$
\rho(s)=\frac{e^{2(n-1)}}{s^{2(n-1)}}=c(n)s^{-2}\,.
$$


\section{\bf Asymptotics for weighted nonlinear diffusion}
\label{sect.wpme}

The study of radial solutions of the PME in hyperbolic space is thus reduced to the study of the weighted nonlinear diffusion equation (\ref{eq.wpme}) with a special type of weight that has a regular behaviour at $s=0$ but degenerates as $s\to\infty$ with an inverse quadratic rate. Complicated non self-similar expressions are to be expected from the analysis of the transformed problem in view of the already complicated expressions for the fundamental solution of the linear heat equation in $\Hn$.

The mathematical study of equations of the form
\begin{equation}\label{eq.5.1}
\rho(x)\,\partial_t {\widehat u}=\Delta_x {\widehat u}^m\,, \quad m>1,
 \end{equation}
 was initiated by Kamin and Rosenau in two seminal papers \cite{KRos81, KRos82}, and continued by a number of authors; no assumption of radial symmetry was used in those studies. The theory is usually known as nonlinear diffusion in inhomogeneous media and $\rho(x)$ is called the weight and represents the mass density of the medium. The previous authors also showed that in the presence of highly decaying weights, like $\rho(x)\le C|x|^{-\gamma}$ with $\gamma>n$, the long-term behaviour of the solutions completely departs  from the non-weighted case; indeed, solutions tend to stabilize to a constant value all over the space. This and other surprising results have led to a keen interest in the subject with contributions by many authors, like \cite{Eidus90, GHP97, GKKV03, KK93, RV2006a, GMP2013, GMPu2013, Pu12}.

The properties and long-time behaviour of the solutions in the presence of weights that decay in a slower way  than $|x|^{-n}$ have been investigated in particular by Reyes and the author in a series of papers \cite{RV2006, RV2008, RV2009}. Thus, for weight functions that degenerate mildly infinity, like $\rho(s)\sim s^{-\gamma}$ with $0\le \gamma<2$, we consider the long-term behaviour of solutions with
nonnegative initial data in the natural physical space
\begin{equation}\label{wpm3.L1rho}
{\widehat u}(x,0)={\widehat u}_0(x)\in L^1_{\rho}(\re^N)
\end{equation}
where $ L^1_\rho(\re^N)= L^1_\rho(x)(\re^N,\rho(x)dx)$. This behaviour takes the form of convergence to self-similar functions of modified Barenblatt type, as established in \cite{RV2009}. The penetration rate for long times is power-like, $L(t)\sim t^{\mu}$ with precise exponent $\mu=1/(n(m-1)+2-m\gamma)$.

 The so-called intermediate case $\gamma\in (2,n)$ was then studied by Kamin and  the same authors in \cite{KRV2010}. Taking a similar class of data, the asymptotic behaviour is now different and more complex, in the form of matched asymptotics. It is first  shown that nontrivial solutions to the problem have a long-time universal behavior in separate variables of the form ${\widehat u}(x, t)=t^{-1/(m-1)}W(x),$ where $V = W^m $ is the unique bounded, positive solution of the sublinear elliptic equation $-\Delta V = c \rho(x)V^{1/m}$ in $\ren$ vanishing as $|x|\to\infty$.
Such a long-time behavior of ${\widehat u}$ is typical of Dirichlet problems on bounded domains
with zero boundary data. It strongly departs from the behavior in the case of
slowly decaying densities. An additional observation it that a separate-variable solution cannot explain the movement of the free boundary. Actually, the ``outer behaviour'' needs another self-similar construction provided also in \cite{KRV2010} which gives a penetration rate $L(t)\sim t^{\mu}$ as before when $\gamma $ is near 2, but the outer description fails for larger $\gamma$.

\subsection{\bf The borderline case}

\noindent The kind of weights that we are interested in the  present application to the hyperbolic context is the borderline case where the weight decays at long distances like $c|x|^{-2}$. A change of character of the equations occurs at this exponent, so the asymptotic description is rather special, and gives the results we are looking for. This change of character is well-known also in the theory of linear heat equations with weights because it corresponds to dimensional considerations (the dimensions of the weight match the dimensions of the Laplacian operator). A key feature is that it involves logarithmic scales.

In the paper \cite{NR2013} my former students Nieto and Reyes studied the long-time behavior of non-negative, finite-energy solutions to the initial value problem for the Porous Medium Equation with variable density, i.e. solutions of equation (\ref{eq.5.1}) with $m>1$,  $n\ge 3$ and nonnegative initial data ${\widehat u}(x,0)={\widehat u}_0(x)\in L^1_{\rho}(\re^N)$. The authors  assume the following set of hypotheses (H) on the weight

\ $\rho\in  C^1(\ren) \cap L^\infty(\ren); \rho > 0,$

\ $\rho(x)|x|^2 \to  1 $ in  $\ren$ as $|x|\to\infty$,

\ $|\nabla \rho(x)|\le C\, |x|^{-3}$ for some  $C> 0$,

\noindent and then show that the long time  behavior of the solutions can be described in terms of a family of source-type solutions of the associated singular equation \ $|x|^{-2}{\widehat u}_t=\Delta_x({\widehat u}^m)$. Indeed, the singular problem  admits the following family of ``Barenblatt solutions'' with a logarithmic singularity at the origin
\begin{equation}
{\widehat U}_A(x, t) = \left[
\frac{\log(A\,t^{\beta}/|x|)}{m(n-2)t} \right]_+^{1/(m-1)}
\end{equation}
where
$$
\beta =\frac1{(n-2)(m-1)}\,.
$$
There is a one-to-one correspondence between the initial weighted norm (usually called energy in this context)
$ E =\int_\ren U(x)|x|^{-2}dx$,  and the free constant $A >0$, of the form
$$
E=k(n,m)A^{n-2}.
 $$
Let us remind the reader that ${\widehat U}_E$ is not a solution of the original equation, but it solves the limit equation with singular weight. Recall also that all solutions are obtained from the solution with $A=1$ by the rescaling $U_A(x,t)= B\,U_1(x,B^{m-1}t)$, with $B=A^{n-2}$. Here is the result of  \cite{NR2013} that we need (after some simplifications of the notation).

\begin{thm} \label{mainteo} Let  $n\ge 3$ and assume that assumptions {\rm (H)} on $\rho$ hold. Let  ${\widehat u}_0 \in L^1_\rho(\ren), \quad {\widehat u}_0 \ge 0$ with $E= ||u_{0}||_{L_{\rho}^{1}}$, and let ${\widehat u}$ be the unique weak solution to the Cauchy problem for Problem {\rm (\ref{eq.5.1})-(\ref{wpm3.L1rho})}. Let ${\widehat U}_{A}$ the self-similar solution with  limit energy $E=k(n,m)A^{n-2}. $ Then,
\begin{equation} \label{est.R1}
\lim _{t\rightarrow \infty} t^{\alpha(p)}||{\widehat u}(\cdot ,t) -
{\widehat U}_{A}(\cdot, t)||_{L^{p}_{\rho}(\mathbb{R}^{n})} = 0
\end{equation}
 for each $p$ with $1 \leq p < \infty$. Here $\alpha(p)= (p-1)/(p(m-1))$.
When $p=\infty$ such convergence fails near $x=0$ and we have for
any $d>0$,
\begin{equation} \label{est.R2}
\lim _{t\rightarrow \infty } t^{1/(m-1)}||{\widehat u}(\cdot ,t) -
{\widehat U}_{A}(\cdot, t)||_{{L^{\infty}(\{|x|\geq dt^{\beta}\})} } = 0.
\end{equation}
\end{thm}

\medskip

\noindent $\bullet$ {\bf New results}

\noindent In our application below we will also need to estimate the location of the free boundary in an accurate way. Therefore, we prove next the following result.

\begin{thm} \label{support.wpme} Let us assume the situation and notations of the previous theorem, and let us consider a solution with  compactly supported initial data and initial energy $E$. Then for large enough times the support of the solution is a star-shaped domain around the origin with maximal and minimal radii ${\widehat R}_-(t)$, ${\widehat R}_+(t)$, such that
\begin{equation}
\lim_{t\to\infty}{\widehat R}_\pm(t)t^{-\beta}=  A.
\end{equation}
In particular, the shape of the support converges to a ball.
\end{thm}

\noindent {\sl Proof.} (i) The lower bound is easy. Due to the uniform convergence (\ref{est.R2}) for all large $t$ we have $u(x,t)>0$ if
$$
\log(A\,t^{-\beta}/|x|) \ge \ve\,.
$$
This means  that $|x|\,t^{\beta} \le Ae^{-\ve}$ implies that $u(x,t)>0$ if $t$ is large enough.

\noindent (ii) The upper bound needs more work. We have to eliminate points in the support such that
$\xi=x\,t^{-\beta}$ is much larger than $A$ for large $t$. Assuming for contradiction that such points exist for a sequence of  times going to infinity, we use a supersolution construction to eliminate them by showing that the excess over $A$ must be reduced in time.

Take any such point where the maximum distance of the support reaches the boundary of the ball of radius $A_1t_0^{\beta}$ with $t_0$ large enough and $A_1\ge k A$, $k>1$. We know by the convergence theorem that
$$
{\widehat u}(x,t)\le \ve t^{-\alpha} \quad \mbox{ for }  \ |x|\ge A\,t^{\beta}, \  t\ge t_0\,,
$$
and also that ${\widehat u}(x,t_0)= 0$ for $|x|\ge A_1\,t_0^{\beta}$. We consider the exterior cylindrical domain $\Omega=\{(x,t): |x|\ge A_1\,t_0^{\beta}, \ t\in (t_0,t_0+h)\}$ and use as a supersolution a function of the form
$$
{\widehat U}_{A'}(x,t+\tau)=
\left(\frac{\log(A'(t+\tau)^{\beta}/|x|)}{m(n-2)(t+\tau)}\right)_+^{1/(m-1)}
$$
with $\tau$ very large and $A'$ to be chosen in order to compare with our solution on the parabolic boundary of $\Omega$. Such comparison is needed at the lateral boundary where $|x|=A_1\,t_0^{\beta}$. If the support of  ${\widehat U}_{A}(x,t)$ does not reach this lateral boundary a sufficient condition  is ${\widehat U}_{A'}\ge \ve t^{-\alpha}$. This situation happens if $A_1\ge kA >A$ and $h/t_0$ is small (depending only on the quotient $A_1/A$). Therefore, we need:
$$
\log(A'(t+\tau)^{\beta})- \log(A_1t_0^\beta)\ge \ve^{m-1} (t+\tau)/t
$$
for $t_0\le t\le t_0+h$. We take
$$
\log A'= \log A_1+ \beta\log (t_0/(t_0+\tau))+ c \ve^{m-1}.
$$
We will assume that $\tau$ is like a multiple of $t_0$.   In that case, by comparison with ${\widehat U}_{A'}$ we get an interface estimate for $t_0\le t\le t_0+h$ of the form
$$
{\widehat R}_+(t)t^{-\beta} \le A'(t+\tau)^\beta t^{-\beta}=A'(1+\tau/t)^\beta.
$$
Using the chosen value of $A'$ we get
$$
\log({\widehat R}_+(t)t^{-\beta})\le \log A_1 + \beta\log (t_0/(t_0+\tau))+ c\ve^{m-1} +\log((t+\tau)^\beta t^{-\beta}).
$$
so that if we call $\xi_+(t)=\log({\widehat R}_+(t)t^{-\beta})$ we have
$$
\xi_+(t)-\log (A_1) \le \beta\log(1+(\tau/t)) -\beta\log (1+(\tau/t_0)) + c\ve^{m-1}.
$$
and we want to prove that $\xi_+(t)$ has strictly negative increments for convenient $h>0$. We will take $h$ as  a small fraction of $t_0$ but always much larger than $\ve^{m-1}$ to kill the last error term. In that case we get approximately
$$
\xi_+(t_0+h)-\log (A_1) \le -\frac{\beta\tau}{t_0t}h+c\ve^{m-1}\le -Ch.
$$
This proves the result that $\limsup_{t\to\infty}\xi_+(t)\le \log (A)$. \qed

Finally, we need an estimate about what happens for $x\sim 0$, where we know that the solution is bounded at every time, but the model solution $ {\widehat U}_A$ has an asymptote. This is usually called the inner analysis, or analysis near the singularity of the model.

\begin{thm} \label{upper} Under the same conditions, there is a constant $C>0$ such that
\begin{equation}
| {\widehat u}(x,t)|^{m-1}\le C\frac{\log t}{t}\,,
\end{equation}
for all large $t>0$ and all $x\in\ren$.
\end{thm}

\noindent {\sl Proof}. We start from estimate ${\widehat u}_t\ge -{\widehat u}/(m-1)t$, see (\ref{beni.est}), which is based on homogeneity and is valid for the equation for ${\widehat u}(x,t)$ as well as the equation for $u(r,t)$. Using this estimate and the equation  we arrive at the inequality
$$
-\Delta {\widehat u}^m \le \frac{\rho}{(m-1)t}{\widehat u}\,,
$$
which can be used to obtain an elliptic estimate as follows. Indeed, the right-hand side is less than
$(\|{\widehat u}(\cdot,t)\|_\infty/(m-1)t)\rho(x)$. Since we know that
$$
\rho(x)\le C\max\{1, r^{-2}\}
$$
We can then estimate ${\widehat u}^m$ by solving the equation
 $$
 -\Delta \Phi=\max\{1, r^{-2}\} \quad \mbox{in \ } B_R(0)\,,
 $$
with $\Phi(x)=0$ for $|x|=R$. We easily get a radial $\Phi(s)$ which is decreasing in $s=|x|$
and such that
$$
-\Phi'(s)\le C\max\{r, r^{-1}\}
$$
and then $\Phi(0)\le C\log R$ if $R\gg 1$.

Choosing $t$ large and applying comparison between ${\widehat u}^m(x,t)$ and $\Phi(x)$ in the ball $ B_R(0)$, where $R={\widehat R}_+(t)$ is the maximum radius of the free boundary that we have just calculated, we get
$$
{\widehat u}^m(x,t)\le \frac{C}{(m-1)t} \Phi(x)\le C \frac{\|{\widehat u(\cdot,t)}\|_\infty}{(m-1)t}(\beta \log(t)+\log A)\,,
$$
which implies the result. \qed

This  estimate must be optimal in view of the equivalence with the hyperbolic problems and our estimate (\ref{subsupest}). This is convenient because obtaining the lower bound in a direct way is not easy.

\section{Application to hyperbolic space}
\label{sect.app}

We return to hyperbolic space and  take a radial initial datum $u_0(r)\in L^1(\Hn)$, $u_0(r)\ge 0$, for $n\ge 3$, and transform it by the change of variables of Section \ref{sect.change} into a function of $s$, ${\widehat u}_0(s)$. We have
$$
M=\int_{\Hn} u_0(r)\,d\mu(x)=\iint u_0(r)(\sinh r)^{n-1}dr\,d\omega_{n-1}
=\iint {\widehat u}_0(s)\rho(s)s^{n-1}ds\,d\omega_{n-1}\,,
$$
so that ${\widehat u}_0\in L^1_\rho(\ren)$. This constant is what we called $E$ in the nonlinear diffusion study of the previous section.  We pass from solutions ${\widehat u}(s,t)$ of the weighted PME equation (\ref{eq.5.1}) to solutions $u(r,t)$ of equation (\ref{eq.hpme}) by means of the transformation (\ref{rtos}). In order to apply the  results of the previous section to find the asymptotic behaviour of ${\widehat u}(s,t)$,  we need to check the necessary (H) conditions in the problem coming from the transformation from the hyperbolic space. Only the last condition needs careful checking: we have seen that as $s\to\infty$
$$
\rho(s)=(\sinh r)^{2(n-1)}s^{-2(n-1)}\sim c_1(n)s^{-2}\,,
$$
 and also $ds/dr\sim c_1 s$ with $c_1=(n-1)/(n-2)$, so that as $s\to\infty$ we get
$$
\frac{d\rho(s)}{ds}=2(n-1)\frac{(\sinh r)^{2(n-1)-1}\cosh r}{s^{2(n-1)}}\,\frac{dr}{ds}
- 2(n-1)\frac{(\sinh r)^{2(n-1)}}{s^{2(n-1)+1}}\,.
$$
$$
=2(n-1)\frac{(\sinh r)^{2(n-1)}}{s^{2(n-1)+1}}(\frac{\coth r}{s'(r)}-1)\sim -\frac{2\rho(s)}{s}
$$
as expected. On a more technical level, in (H) we have imposed the normalization condition
$\rho(s)s^2\to 1$ as $s\to\infty$, which does not happen according to (\ref{rho.def}) because of a factor $c_1(n)$, that can be absorbed into the time variable, using new time $t'=t/c_1(n)$.

As a consequence, we  can  describe the long-time behaviour of radial solutions of the porous medium equation in hyperbolic space in terms of the explicit function ${\widehat U}_A(s,t')$, that after the transformation reads
\begin{equation}
{\widetilde U}(s,t')=\left( \frac{\log(A(t')^{\beta}/s(r))}{m(n-2)t'}\right)_+^{1/(m-1)}
\end{equation}
It is quite interesting to notice that this is just the approximate supersolution $\widetilde U$ that we have constructed in (\ref{approx.sol}) after the formal substitution $s=e^{(n-1)r/(n-2)}$ and adaptation of the constants. It is not immediately clear that the constant are  just the ones of Theorem (\ref{thm.main}), but that is so by after checking that both functions $\widetilde U$ and $\widetilde U$ agree.

Since the asymptotic result concerns  large values of $r$ and $t$, we get as a consequence the following result, valid for in particular the fundamental solutions.

\begin{thm}\label{teo61}  Under the stated assumptions on the radial solutions of the HPME,  we have as $t\to\infty$ and for $r\ne 0$
\begin{equation}\label{conv.hpme}
t\,u(r,t)^{m-1}\sim a\left(\,\gamma\log t - r + \gamma(m-1)\log (A)+c\right)_+\,,
\end{equation}
with $\gamma$ and $a$ as in Theorem   {\rm \ref{thm.main}} and uniform convergence on sets of the form $r\ge \gamma\log(t)+d$. The convergence of $u(\cdot,t)$ towards $\widetilde U(r,t;A)$ takes place in $L^p(\Hn)$ for all $p\in[1,\infty)$, but not in $L^\infty(\Hn)$ of course.
\end{thm}

 We need to supplement this bulk information with a precise information of what happens near $s=\infty$ and near $s=0$.

\begin{thm} \label{teofb} Under the above conditions and notations, if moreover the initial data have compact support then for large time the support looks like a ball of radius $R(t)$ and
\begin{equation}\label{eq.fb.lim}
\lim_{t\to\infty}\frac{R(t)}{\log (t)}=\gamma.
\end{equation}
Moreover,
\begin{equation}
\lim_{t\to\infty}e^{R(t)}\,t^{-\gamma}=c_3(n,m)A^{(n-2)/(n-1)}.
\end{equation}
\end{thm}

\begin{thm} \label{teo.upper} Under the same conditions, there is a constant $C>0$ such that
\begin{equation}
| u(r,t)|\le C\left(\frac{\log t}{t}\right)^{1/(m-1)}\,,
\end{equation}
for all large $t>0$ and all $r>0$.
\end{thm}

These results cover all the statements of our main result, Theorem \ref{thm.main}, extending them in different ways.

\subsection{The behaviour of more general solutions}

We have seen in this section how the results we were looking for fundamental solutions extend to
the class of nonnegative radial solutions with compactly supported data. In the case where the solutions are not radial but still compactly supported, the standard maximum principle argument allows to sandwich such a solution between two radial solutions with same origin point $O$, and then the two first theorems hold the weaker version of bounds from above and below with respect to expressions with different energies, i.\,e., with a large $A_1$ from above and a small $A_2$ from below. In Theorem \ref{teofb} the first formula, (\ref{eq.fb.lim}), holds literally, but the second needs constants. In the last theorem we may use two constants $C_1$ and $C_2$.

\section{Analysis in two space dimensions}\label{sect.2d}

An upper bound for the solutions follows from the construction of the approximate solution in Section \ref{sect.aaas}. This argument is the same as in higher dimensions. We get
$$
U(r,t)\le t^{-1/(m-1)}(a\gamma \log(t)+b)^{1/(m-1)}
$$
with $a=1/m(m-1)$, $\gamma=1/(m-1)$ and  $b=b(m,n)$. It follows that
$$
R(t)\le \frac1{m-1} \log(t)+ b.
$$
According to what we have seen for $n\ge 3$ these bounds should be sharp.

\medskip

\noindent $\bullet$ However, the question of sharp estimates from below in 2D is more difficult because the equivalence and transfer method of previous section does not work. Let us see why: revisiting what was said in Section \ref{sect.change}, we find that the change of variables is now  given by
$$
\frac{ds}{dr}=\frac{s}{\sinh r}
$$
This means that for $r\sim 0$ we have $s(r)\sim r$ while for $r\gg 1$ we get
$$
\log (s(r))\sim c\,e^{r}.
$$
After some computations as before, we get for $\widehat u$ the same equation:
$\rho(s)\,\partial_t {\widehat u}=\Delta_s {\widehat u}^m,$
but now the weight is
$$
\rho(s)=(\sinh r)^{2}/s^{2}\,,
$$
so that  $\rho(s)\sim 1$ for $s\sim 0$ while $\rho(s)\sim (\log\log (s))/s^2$ for $s\gg 1$.
Contrary to the case $n\ge3$ we are not dealing with a small perturbation of the inverse square weight; in this case the correction is a $log$-$log$ term.

\noindent $\bullet$ There some lower estimates that can be obtained by remarking that radial solutions for $n=3$ are subsolutions for $n=2$, and this implies bounds from below where only the value of the coefficients changes. Thus, we get
$$
R(t)\ge \frac1{2(m-1)} \log(t)+ b'.
$$
But in fact, we can run the equivalence method of Sections \ref{sect.change} and \ref{sect.wpme} for radial solutions in any non-integer dimension (we refrain here for lack of space from entering into the detailed justification of this assertion, on the confidence that the reader will check it if needed). In this way we obtain lower bounds for the minimal radius of the support of the form $\lim_{t\to\infty}R(t)/\log (t)\ge \gamma_n $ with $\gamma_n=\/(m-1)(n-1)$. Passing to the limit $n\to 2$ and combining with the upper estimate we get

\begin{thm} Let $u$ be a solution of the HPME in $\mathbb{H}^2$ with nonnegative, bounded and compactly supported data. Then for all large $t$ the support of the solution looks lie a ball with radius $R(t)$ and
\begin{equation}
\lim_{t\to\infty}\frac{R(t)}{\log (t)}=\frac1{m-1}\,.
\end{equation}
From this we conclude like in previous sections that
\begin{equation}
\|u(\cdot,t)\|_\infty^{m-1} \sim \frac{\log t}{t}\,.
\end{equation}
\end{thm}

From what we have said above we can also deduce by comparison a weaker version of the estimate
$$
t\,U(r,t)^{m-1}\sim  a \left(\gamma \log t- r \right)_+
$$
with $a=1/m$ and $\gamma=\/(m-1)$,, but not the strong convergence form stated in Theorem \ref{thm.main} for $n\ge 3$.

\section{Extension. The evolution p-Laplacian equation}\label{sect.pl}

We  consider  in this section the nonlinear evolution equation called the $p$-Laplacian equation
\begin{equation}
\partial_t u=\nabla\cdot(|\nabla u|^{p-2}\nabla u)
\end{equation}
with $p>2$.  It is well known that for this equation  there exists in the Euclidean case a family of Barenblatt solutions with self-similar form, parallel to the ones of the PME case, \cite{Bar1952}, and they explain the asymptotic behaviour \cite{KV88}. Much of the theory of the equation, that can be consulted in \cite{DB93}, bears great similarities with the porous medium equation, and that parallelism is stressed for instance in \cite{VazSm2006}, Chapter 11.

This equation has been less studied on Riemannian  manifolds. A basic analysis of the propagation properties of the $p$-Laplacian equation in the degenerate case $p>2$ posed on a complete Riemannian  manifold $(M,g)$ is performed in \cite{Dek05} where it is proved that non-negative bounded solutions of this equation have finite propagation speed, see also \cite{Dek08}\footnote{The author thanks Fabio Punzo for providing this information.}.

We are going to see that our ideas for the PME extend to the $p$-Laplacian equation posed in the hyperbolic space in the parameter range $p>2$ that corresponds to slow propagation and existence of free boundaries. We will be interested in radial, nonnegative and compactly supported solutions. The $p$-Laplacian equation in the hyperbolic space (HPLE) can be written in  radial coordinates as
\begin{equation}
\begin{aligned}
u_t &=(\sinh r)^{1-n}\left((\sinh r)^{n-1}|u_r|^{p-2}u_r\right)_r\\
& =(|u_r|^{p-2}u_r)_r+(n-1)\coth(r)|u_r|^{p-2}u_r\,.
\end{aligned}
\end{equation}
Note that in the limit $p=2$ we get again the heat equation.

As a key calculation to start the programme that we have just developed  for the HPME, we present here an approximate solution for the HPLE. The approximate solution is compactly supported, but since it will be spread over a large part of the space, we may again replace for all practical purposes  the coefficient $\coth(r)$ by 1, and consider exact solutions of the following approximate equation
\begin{equation}
\partial_t u=(|u_r|^{p-2}u_r)_r+(n-1)|u_r|^{p-2}u_r\,.
\end{equation}
We want to try a solution of this equation of the form
\begin{equation}
{\widetilde U}(r,t)= t^{-1/(p-2)}F(\gamma \log t-r+b)\,.
\end{equation}
Then $F(\xi)$ must satisfy
\begin{equation}\label{profileFp}
-\frac1{p-2}F(\xi)+\gamma F'(\xi)= (|F'|^{p-2}F')'-(n-1)|F'|^{p-2}F'\,,
\end{equation}
where primes denote derivative with respect to $\xi\in\re$. In this way the time dependence is eliminated. We  now propose as  profile function  \
$F=a(\xi)_ +^{(p-1)/(p-2)}$. The coefficients $a$ and $\gamma$ are to be determined. Inserting it into   (\ref{profileFp}) we get the conditions, to be satisfied for $\xi\in (-b/a,\gamma \log(t))$,
$$
-\frac{a}{p-2}\xi^{(p-1)/(p-2)} + \frac{\gamma a(p-1)}{p-2}\xi^{1/(p-2)}=
$$
$$-(n-1)\left(\frac{a(p-1)}{p-2}\right)^{p-1}\xi^{(p-1)/(p-2)} + a^{p-1}\left(\frac{(p-1)}{p-2}\right)^{p}\xi^{1/(p-2)}\,.
$$
After simplifying
$$
-a\xi +\gamma a(p-1)= -(n-1)(p-2)a^{p}\left(\frac{p-1}{p-2}\right)^{p-1}\xi+
(p-2)a^{p-1}\left(\frac{a(p-1)}{p-2}\right)^{p}.
$$
This is satisfied if
$$
a^{p-2}=  \frac{(p-2)^{p-2}}{(n-1)(p-1)^{p-1}}\,, \quad \gamma=\frac{1}{(p-2)(n-1)}\,.
$$
We see that  $b\in \re$ is a free parameter (corresponding to space translation invariance; this constant can be used to adjust the mass of the solution at any given time). Summing up,  the approximate solution is
\begin{equation}
{\widetilde U}(r,t)= at^{-1/(p-2)}(\gamma \log t-r+b)^{(p-1)/(p-2)}\,.
\end{equation}

\noindent $\bullet$ From this and convenient theory that must be developed, we must conclude the typical long-time behaviour of the free boundaries \begin{equation}
R(t)=\gamma \log(t) + l.o.t.
\end{equation}
($ l.o.t.$ meaning lower order terms), as well as the sup estimate
\begin{equation}
\|u(\cdot,t)\|_\infty\le a{t^{-1/(p-2)}}\,\log (t)^{(p-1)/(p-2)}+ l.o.t.
\end{equation}
We invite the reader to continue this study, that can be extended to the the doubly degenerate equation
$$
\partial_t u=\nabla\cdot(|u|^{m-1}|\nabla u|^{p-2}\nabla u)
$$
with $p>2$, $m>1$. See \cite{SV13} for some recent work for this equation on Euclidean space.


\section{Comments, extensions and open problems}\label{sect.comm}

\noindent $\bullet$ Our results demonstrate a main feature of porous medium flow in hyperbolic space, the very slow propagation at long distances. Thus, the growth rate of the support for large times is qualitatively smaller than in the Euclidean case, and depends on dimension only in a secondary way, through the coefficient. The same applies to the decay of the sup norm of the solution. On the other hand, the size of solution and expansion of support are larger than in the Dirichlet problem posed in a ball $B_R(O)$. The hyperbolic space rates are a kind of interpolation between both situations, surprisingly closer to the Dirichlet problem than to the Cauchy problem in the whole space $\ren$.

\noindent $\bullet$ For the  sake of completeness and comparison, it is interesting to consider the porous medium equation posed on an $n$-dimensional sphere $\Sn$ in order to complete the picture with a manifold of constant positive curvature. If $f(r)$ is a radial function (with $r$ the geodesic distance  taken from a point $N$ called the North Pole; in common language, the latitude) the radial Laplacian is
\begin{equation}
\Delta_{\Sn} f=(\sin r)^{1-n}\left((\sin r)^{n-1}f'(r)\right)'\,.
\end{equation}
It is easy to construct a fundamental solution for the PME on such manifold, to prove that it is radial, and also that is has compact support for some time $0<t<T$ until the solution spreads to the whole sphere. After that time the fundamental solution proceeds to stabilize towards a constant. In the context of the  diffusion theory of Section \ref{sect.wpme} in inhomogeneous media, this would correspond to the stabilization results by Kamin and Rosenau \cite{KRos81, KRos82} and the disappearing interfaces studied in \cite{GHP97}.

\noindent $\bullet$ In our study we have considered the PME exponents $m>1$. For $m<1$ we get the Fast Diffusion Equation that has quite different properties when posed in Euclidean space, and also in manifolds   \cite{BGV08}.  In particular, nonnegative solutions are indeed positive and exhibit no free boundaries.  Also, a large class of solutions extinguish in finite time.

In a recent paper, Grillo and Muratori \cite{GM2013} have studied the FDE  posed on the hyperbolic space $\Hn$. They also treat radial nonnegative solutions and they need the  restriction $m\in(m_s,1)$, $m_s=(n-2)/(n+2)$, $n\ge 2$. The authors  establish that  near the extinction time $T>0$ the behaviour of the solution is expressed in terms of a separable solution of the form
$$
U(r,t)=(1-t/T)^{1/(1-m)}V^{1/m}(r),
$$
where $V$ is the unique positive and radial energy solution to the elliptic equation $-\Delta_g V=cV^{1/m}$ for a suitable $c>0$. Such semilinear elliptic equations were thoroughly studied and classified in \cite{BGGV13}, \cite{MS2008}. In particular, the solution of the fast diffusion equation is bounded above and below near the extinction time $T$ by multiples of $(1-t/T)^{1/(1-m)}e^{-(N-1)r/m}.$ This is to be compared  to the quite different results of the linear heat equation, which appears in the limit $m\to 1$.

\noindent $\bullet$ A standard TW is an expression of the form $U=F(x_i-t)$, while by Generalized Traveling Waves we mean expressions of the form $U=a(t)F(b(t)(x_i-c(t)))$ with time dependent functions $a(t)$, $b(t)$ and $c(t)$ to be determined. If $a(t)$ is decreasing we say that the wave is damped; $c'(t)$ is the variable speed. Our explicit solutions (\ref{gtw.22}),  (\ref{gtw.mn}), conform to that pattern only if logarithmic variables are introduced (which is natural in view of our results). TW solutions in such generality were considered for the PME theory in the Euclidean space, cf. \cite{Vbook}, subsection 4.7, but they have not played a big role. Other options exist in the literature for solutions which have some essential properties of TW's: thus, Generalized Transition Waves have been discussed in \cite{BH12},  where  further references can be found.

\noindent $\bullet$ Here is an alternative way to find the generalized TW (\ref{gtw.22}). We want to find a solution with internal speed $1/t$ everywhere (on the positivity set) and depending only on the variable $y$ in the Poincar\'e half plane representation. This determines the solution. Indeed, the speed condition means that $P_s=-1/t$, i.\,e., $P_y=1/(ty)$.
Using the equation and the fact that $P=2U$ we have
$$
P_t=2U_t=4y^2(UU_y)_y= y^2\left(\frac{P}{ty}\right)_y=
 y^2\frac{P_y}{ty}- y^2\frac{P}{ty^2}
$$
Therefore, $tP_t+P=1/t$, so that $tP=\log(t)+f(y)$. Since we also have $tP_y=1/t$ we get the final form $tP=\log(cty)$ for all points where $P>0$. The author finds this derivation very mechanical and insightful. A similar derivation could be done for formula (\ref{gtw.mn}) with general $m$ and $n$.

\noindent $\bullet$  The traveling wave behaviour for the renormalized equation (\ref{renorm1b}) looks promising and could lead to better understanding. Note that is corresponds to propagation with logarithmic rate in the original time variable. A similar situation was found by the author in the study of propagation of the PME in tubular domains in the Euclidean space (with zero Dirichlet conditions), see \cite{VazTube}. In that case the TW is bounded, there is no logarithmic growth in time.

\noindent $\bullet$  Our level of detail in the theory of the PME in hyperbolic space has been rather sketchy, but we feel that it is sufficient, since we mention some basic references.   A longer treatment would correspond to a more expository text that should be written, and  is not justified here, where the goal is the construction of special solutions and solving the question of precise long-time behaviour.

\noindent $\bullet$ As for other questions worth studying, let us mention:

(i) We have not discussed in this paper the detailed  behaviour of the solution in the inner region, i.\,e., for $e^r \le \varepsilon  t^{\gamma}$. In particular, the important question of obtaining the sharp estimate of $\|U(\cdot,t)\|_\infty$ as $t\to\infty$ has not been settled. The expected way of doing the analysis of this inner behaviour is through matched asymptotics, taking as example the analysis of the paper \cite{DGLV1998} for the heat equation with exponential reaction. The adaptation to the present situation is not easy, so it is not addressed  here.

(ii) We leave untouched the questions of precise geometry and regularity of the free boundaries of more general solutions.

(iii) The sharp estimates for the HPME posed in 2 dimensions are an interesting open problem too.

(iv) We have studied the existence, uniqueness and we have obtained sharp estimates on the behaviour of the fundamental solutions and other solutions in the context of radial solutions in hyperbolic space. It will be interesting to get similarly sharp results for non-radial solutions.

(v) It would be interesting to know if there are more explicit formulas for the fundamental solutions than the ones provided in Theorem \ref{thm.main}.

(vi) This is a first step into the study of finite propagation of porous medium type on manifolds with negative curvature. This should be pursued with the aim of better understanding the influence of geometry on finite propagation.

(vii) The existence of the special solution \eqref{gtw.22}--\eqref{gtw.mn}, studied in Section \ref{sect.tw}, poses the problem of finding the optimal class of uniqueness for nonnegative solutions of the Cauchy Problem for equation \eqref{hpme} in hyperbolic space. Such question has been extensively studied in $\ren$: in the heat equation case er refer to Widder \cite{Wid}; in the PME case cf. the fundamental work by B\'enilan et al. \cite{BCP84} or the general presentation in the monograph \cite{Vbook}.  The same problem is posed for the $p$-Laplacian counterpart.


\section*{Appendix}

\noindent $\bullet$ {\sl Calculation for the generalized traveling wave of Section \ref{sect.tw}.} We want to check that the explicit function  given by
\begin{equation}
U(x,y,t)^{m-1}=\frac{(a\log(ct^{\gamma}y))_+}{t}
\end{equation}
is a solution of equation $\partial_t u=\Delta_{\Hn} (u^m)$ for all  $m>1$, $n>1$, when  $a=1/m(n-1)$ and $\gamma=1/(m-1)(n-1)$.

\noindent {\sl Proof.} Let $\Omega=(a\log(ct^{\gamma}y))_+$ and $\alpha=1/(m-1)$. Then when the solution is positive we have
$$
U=t^{-\alpha}\Omega^{\alpha}, \qquad
U_t=-\alpha t^{-\alpha-1}\Omega^{\alpha}+a\gamma t^{-\alpha-1}\Omega^{\alpha-1}
$$
$$
U^m=t^{-m\alpha}\Omega^{m\alpha}, \qquad
(U^m)_y= m\alpha t^{-m\alpha}\Omega^{m\alpha-1}\frac{a}{y}=m\alpha t^{-m\alpha}\Omega^{\alpha}\frac{a}{y}
$$
$$
(U^m)_{yy}=-m\alpha t^{-m\alpha}\Omega^{\alpha}\frac{a}{y^2}+ m\alpha^2 t^{-m\alpha}\Omega^{\alpha-1}\frac{a^2}{y^2}
$$
so that the equation translates into
$$
y^2(U^m)_{yy}-(n-2)(U^m)_y=-m\alpha a t^{-m\alpha}\Omega^{\alpha}+ m\alpha^2 a^2 t^{-m\alpha}\Omega^{\alpha-1}-(n-2)m\alpha a t^{-m\alpha}\Omega^{\alpha}
$$
$$
= -ma\alpha (n-1)t^{-m\alpha}\Omega^{\alpha}+ ma^2\alpha^2 t^{-m\alpha}\Omega^{\alpha-1}
$$
This is satisfied if $\alpha=ma\alpha (n-1), \qquad a\gamma=ma^2\alpha^2.$
In order to check the FB condition we argue as before.

\medskip

\noindent $\bullet$ {\sl Transformation of the full hyperbolic Laplacian}
\label{sect.nonrad}

There is a legitimate question of trying the transformation into a weighted equation for nonradial functions. As we have said, the full Laplace-Beltrami operator  in hyperbolic  space $\Hn$, $n\ge 2$, can be written as
\begin{equation}\label{eq.hyplap}
\Delta_g \,u(r,\xi) =
\frac1{(\sinh r)^{N-1}}\left((\sinh r)^{N-1}u_r(r)\right)_r+ \frac1{\sinh(r)^{2}}\Delta_{\xi}u(r,\xi)\,.
\end{equation}
where $0<r<\infty$ is the geodesic distance and  $\xi\in {\mathbb S}^{n-1}$. We use again as change rule
\begin{equation}
\frac{ds}{s^{n-1}}=\frac{dr}{(\sinh r)^{n-1}}\,.
\end{equation}
After some easy computations we pass  from the original Laplacian into the form
 $$
\Delta_g \,u(r,\xi) =\frac{s^{(n-1)}}{(\sinh r)^{2(n-1)}}\left(s^{n-1}u_s(s,\xi)\right)_s  +
 \frac1{\sinh(r)^{2}}\,\Delta_{\xi}u(r,\xi)\,.
 $$
In other words, $\Delta_g $ transforms into the differential operator $L$ in the variables $s$ and $    \xi$ given by
 $$
L {\widehat u}(s,\xi) =\frac1{\rho(s)}\left(\frac1{s^{n-1}}(s^{n-1}{\widehat u}_s(s,\xi))_s+ \frac{\mu(s)}{s^{2}} \Delta_{\xi}{\widehat u}(s,\xi)\right)\,,
$$
where ${\widehat u}(s,\xi)=u(r,\xi)$ and
$$
\rho(s)=\left(\frac{\sinh r}{s}\right)^{2(n-1)}\,, \qquad
\mu(s)=\frac{s^2\rho(s)}{(\sinh r)^2}=\left(\frac{\sinh r}{s}\right)^{2(n-2)}\,.
$$
We conclude that the transformation produces a weighted Laplacian in Euclidean space that also distorts the non-radial terms for every $n\ge 3$ since $\mu\ne 1$, but it does not have the distortion effect for $n=2$.

\

\noindent {\large\bf Acknowledgments}.   Work partially funded by Spanish Grant MTM2011-24696. The author would like to thank the Isaac Newton Institute for Mathematical Sciences, Cambridge, for support and hospitality during the Free Boundary programme, Jan. to July 2014, where work on this paper was undertaken.  The text has also benefited from positive suggestions from a number of colleagues.

\vskip 1cm


%

%
\


{\sc Full Address:}

Juan Luis V\'azquez, Departamento de Matem\'{a}ticas,

Universidad Aut\'{o}noma de Madrid, Campus de Cantoblanco,

28049 Madrid, Spain.

E-mail:{{\tt~juanluis.vazquez@uam.es}.

\vskip .5cm


\noindent {\bf Keywords.} Porous Medium Equation, Hyperbolic Space, Fundamental Solutions, Asymptotic Behavior.
\\[.2cm]
{\sc 2010 Mathematics Subject Classification}.
 	35K65,   	
58J35,   	
 	35K55,   	
35B40.   	

\newpage

\end{document}